\documentclass[12pt]{article}
\usepackage[a4paper, total={6in, 9.6in}]{geometry}
\usepackage{amsmath,amssymb}
\usepackage{verbatim}
\usepackage{enumerate}
\usepackage{authblk}

\usepackage{graphicx}
\usepackage{caption}
\usepackage{subcaption}

\newtheorem{theorem}{Theorem}
\newtheorem{corollary}[theorem]{Corollary}
\newtheorem{lemma}[theorem]{Lemma}

\newcommand{\lap}{\mathcal{L}}
\newcommand{\R}{\mathbb R}
\newcommand{\C}{\mathbb C}

\newcommand{\cD}{{}^c\!D}
\renewcommand{\Re}{\mathrm{Re} \;}
\newcommand{\etal}{\mbox{\emph{et al.\ }}}
\newcommand{\hD}{\!D^{\alpha,\gamma}}

\newcommand{\proof}{\noindent{\bfseries Proof. }}
\newcommand{\eproof}{\mbox{\ \rule{.1in}{.1in}}}

\newcommand{\thline}{\rule{5.75in}{.01in}}

\begin{document}
\title{An inverse source problem for a two-parameter anomalous diffusion with local time datum}

\author[1]{Khaled M. Furati}
\author[2]{Olaniyi S. Iyiola}
\author[3]{Kassem Mustapha}

\affil[1,3]{\small King Fahd University of Petroleum \& Minerals, Department of Mathematics \& Statistics,
Dhahran 31261, Saudi Arabia
}
\affil[1]{kmfurati@kfupm.edu.sa}
\affil[3]{kassem@kfupm.edu.sa}

\affil[2]{\small University of Wisconsin-Milwaukee, Department of Mathematical Sciences, WI, USA
}
\affil[2]{osiyiola@uwm.edu}

\maketitle


\begin{abstract}

We determine the space-dependent source term for a two-parameter fractional diffusion problem subject to nonlocal non-self-adjoint boundary conditions and two local time-distinct datum.
A bi-orthogonal pair of bases is used to construct a series representation of the solution and the source term.
The two local time conditions spare us from measuring the fractional integral initial conditions commonly associated with fractional derivatives.
On the other hand, they lead to delicate $2\times 2$ linear systems for the Fourier coefficients of the source term and of the fractional integral of the solution at $t=0$.
The asymptotic behavior and estimates of the generalized Mittag-Leffler function are used to establish the solvability of these linear systems, and to obtain sufficient conditions for the existence of our construction.  Analytical and numerical examples are provided.
\end{abstract}


\thline

\section{Introduction}
Non-integer calculus and fractional differential equations have become an intrinsic tool in modeling different phenomena in many areas such as nanotechnology \cite{baleanu:10b}, control theory of dynamical systems \cite{caponetto:10b,monje:10b}, viscoelasticity \cite{mainardi:10b}, anomalous transport and anomalous diffusion \cite{klags:08b}, random walk \cite{hilfer:95,zhang:06}, financial modeling \cite{scalas:04} and biological modeling \cite{olaniyi:2014}. Some other physical and engineering processes are given in \cite{ortigueira:11b,petras:11b} and more applications can be found in the surveys in \cite{hilfer:00b, kilbas:06b, podlubny:99b}.
In particular, fractional models are increasingly adopted for processes with anomalous diffusion \cite{adams:92,metzler:00,zhou:03}.
The featured role of the fractional derivatives is mainly due to their non-locality nature which is an intrinsic property of many complex systems \cite{hilfer:00}.

In this paper, we consider determining the solution $u$ and the space-dependent source $f$ for the following two-parameter fractional diffusion equation (FDE):
\begin{eqnarray}
\label{eq:heat-equation}
&& \nonumber
\hD u(x,t)-u_{xx}(x,t) = f(x),
\hspace{.4in} x \in (0,1),\quad t \in (0,T), \quad 0< \alpha \leq \gamma \leq 1,
\\ &&
u(x,T_m)=z(x),\quad u(x,T)=h(x),\quad x \in [0,1], \qquad 0 < T_m < T,
\\ && \nonumber
u(1,t)=0,\quad u_x(0,t)=u_x(1,t),\hspace{.4in} t \in (0,T],
\end{eqnarray}
where $z$ and $h$ are square integrable functions. The operator $\hD$  is the generalized Hilfer-type fractional derivative defined by
\begin{equation}
\label{eq:hD-def}
\hD y(t)= D^\alpha \left[ y(t) - I^{1-\gamma} y(0) \,
\frac{t^{\gamma -1}}{\Gamma(\gamma)} \right],
\end{equation}
where
$$
I^\alpha y(t)=\frac{1}{\Gamma(\alpha)}\int_{0}^{t}(t-\tau)^{\alpha-1} y(\tau) \, d\tau,
$$
and
$$
D^\alpha y(t) = D I^{1-\alpha} y(t), \qquad \quad D = \frac{d}{dt},
$$
are the Riemann-Liouville fractional integral and derivative, respectively. Noting that, when $I^{1-\gamma} y \in AC[0,T]$, then $\hD y$ can be written as
\begin{equation}
\label{eq:hD-hilfer}
\hD y(t) = I^{\gamma-\alpha} D^\gamma y(t) = I^{\gamma-\alpha} D I^{1-\gamma} y(t).
\end{equation}
Hence, $\hD$ reduces to the derivative introduced by Hilfer in \cite{hilfer:00evolution} when  $\gamma = \beta(1-\alpha)+\alpha$ with $0<\beta\leq 1$.

Observe that the Riemann-Liouville fractional derivative $D^\alpha$ and the Caputo fractional derivative $\cD^\alpha:= I^{1-\alpha} D$ are special cases of the two-parameter fractional derivative $\hD$ for $\gamma=\alpha$ and $\gamma=1$, respectively. Thus,  $\hD$ is considered as an interpolant between $D^\alpha$ and $\cD^\alpha$.

Furati \etal \cite{furati:14} constructed a series representation of $u$ and $f$ for problem \eqref{eq:heat-equation}, but subject to an integral-type initial condition instead, using a bi-orthogonal system. Unlike in \cite{furati:14}, in the current problem, the value of the solution at some time $t = T_m$ is used rather than the value of the fractional integral of the solution $I^{1-\gamma} u(x,t)$ at $t=0$, which may neither be measurable nor have a physical meaning.
As a result, the construction of the series representation of $u$ and $f$ is not a straightforward extension of the one in \cite{furati:14}.
This is mainly due to the complexity in ensuring the solvability of the arising linear systems and in achieving the necessary lower and upper bounds for showing the convergence of the constructed series.

Inverse source problems for a one-parameter FDE with Caputo derivative have been investigated by many researchers under various initial, boundary and over determination conditions. For a space-dependent source $f$,  
Zhang and Xu \cite{zhang-xu:11} used Duhamel's principle and an extra boundary condition to uniquely determine $f$.
Kirane and Malik \cite{kirane:11} studied first a one-dimensional problem with non-local non-self-adjoint boundary conditions and subject to initial and final conditions. The results were extended to the two-dimensional problem by Kirane \etal \cite{kirane:12}.
\"{O}zkum \etal \cite{ozkum:13} used Adomian decomposition method to construct the source term for a linear FDE with a variable coefficient in the half plane.
In a bounded interval, Wang \etal \cite{wang:13} reconstructed a source for an ill posed time-FDE by Tikhonov regularization method.
A numerical method for reproducing kernel Hilbert space to solve an inverse source problem for a two-dimensional problem is proposed by Wang \etal \cite{wang:13num}.
Wie and Wang \cite{wei:14} proposed a modified version of quasi-boundary value method to determine the source term in a bounded domain from a noisy final data.
Analytic Fredholm theorem and some operator properties are used by Tatar and Ulusoy \cite{tatar:15} to prove the well-posedness of a one-dimensional inverse source problem for a space-time FDE.
Feng and Karimov \cite{feng:15} used eigenfunctions to analyze an inverse source problem for a fractional mixed parabolic hyperbolic equation.
They formulated the problem into an optimization problem and then used semismooth Newton algorithm to solve it. For  a three-dimensional inverse source problem, we refer the reader to the work  by Sakamoto and Yamamoto \cite{sakamoto:11} and  by Ruan \etal \cite{ruan:16}. In relation to above, for the case of  time-dependent source term $f$, see  \cite{aleroev:13,ismailov:16, sakamoto:11ibvp, wei:13, wu:14, yang:15}.

The rest of the paper is organized as follows.
In section \ref{sec:GFML}, we present a brief discussion of the generalized
Mittag-Leffler function and derive some related properties.
In addition, we establish the solvability of a $2\times 2$ linear system with a coefficient matrix involving Mittag-Leffler functions and obtain estimates of these coefficients.
We construct the series representations of the solution and source term in
Section \ref{sec:soln}.
The well-posedness of this construction is proved in Section \ref{sec:results}.
In Section \ref{sec:example}, analytical and computational examples are presented.

\section{Generalized Mittag-Leffler Function}
\label{sec:GFML}

The Prabhakar generalized Mittag-Leffler function \cite{prabhakar:71} is defined as
\begin{equation}
\label{eq:GMLF1}
E_{\alpha,\beta}^\rho (w) = \sum_{k=0}^\infty
\frac{\Gamma(\rho+k)}{\Gamma(\rho) \, \Gamma(\alpha k + \beta)} \, \frac{w^k}{k!},
\qquad w, \beta, \rho \in \C, \quad \Re \alpha > 0.
\end{equation}
Some special cases of this function are the Mittag-Leffler function in one parameter
$E_\alpha(w) = E_{\alpha,1}^1 (w)$, and in two parameters
$E_{\alpha,\beta} (w) =  E_{\alpha,\beta}^1 (w)$.

The function $E_{\alpha,\beta}^\rho$ is an entire function  \cite{prabhakar:71} and therefore bounded in any finite interval.
In addition, for $w,\beta, \lambda \in \C$ and $\Re \alpha > 0$, this function satisfies the following recurrence relations \cite{haubold:11,kilbas:06b}:
\begin{equation}
\label{eq:ML-recurrence}
\lambda w^\alpha \; E_{\alpha,\alpha+\beta}(\lambda w^\alpha) =
E_{\alpha,\beta}(\lambda w^\alpha) - \frac{1}{\Gamma(\beta)},
\end{equation}
and
\begin{equation}
\label{eq:GML-recurrence}
\alpha E^2_{\alpha,\beta}(\lambda w^\alpha) = (1+\alpha - \beta) \;
E_{\alpha,\beta}(\lambda w^\alpha) + E_{\alpha,\beta-1}(\lambda w^\alpha).
\end{equation}
Combining the relations
\eqref{eq:ML-recurrence} and \eqref{eq:GML-recurrence} yields the recurrence relation:
\begin{equation}
\label{eq:E2->E}
\alpha \lambda w^\alpha \; E^2_{\alpha,\beta}(\lambda w^\alpha) =
(1+\alpha - \beta) \; E_{\alpha,\beta -\alpha}(\lambda w^\alpha)
+
E_{\alpha,\beta -\alpha-1}(\lambda w^\alpha),
\end{equation}
for $w,\beta, \lambda \in \C$ and  $\Re \alpha > 0$.


On the real line, we have the following upper bounds.

\begin{lemma}
\label{lem:mlf-bound}
Let $0<\alpha<2$, $\beta \in \R$, and $\rho = 1, 2$.  Then, there is a constant $M = M(\alpha,\beta) > 0$ such that
\begin{equation}
\label{eq:mlf-bound}
\lambda^{\rho} t^{\alpha\rho} \left|E_{\alpha,\beta}^\rho(-\lambda t^\alpha) \right| \leq  M,
\qquad  \quad t \geq 0, \quad \lambda \geq 0.
\end{equation}
\end{lemma}

\proof
When $\rho = 1$, the result is a special case of Theorem 1.6 in \cite{podlubny:99b}.
When $\rho = 2$, the bound follows from \eqref{eq:E2->E}.  \eproof

\bigskip

The function $E^\rho_{\alpha,\beta}$ possesses the following positivity and monotonicity properties \cite{gorenflo:14b}.

\begin{lemma}
\label{lem:mlf-cm}
Let $0 < \alpha \leq 1$, $\lambda > 0$, and $\rho > 0$.
Then the functions
$E^\rho_{\alpha,\beta}(-\lambda t^\alpha)$, $\alpha \rho \leq \beta$, and
$t^{\gamma-1} E^\rho_{\alpha,\gamma}(-\lambda t^\alpha)$, $\alpha \rho \leq \gamma \leq 1$, are positive monotonically decreasing functions of $t>0$.
\end{lemma}
Lemma \ref{lem:mlf-cm} and the recurrence relation \eqref{eq:ML-recurrence} imply the following corollary.
\begin{corollary}
\label{cor:mlf-increasing}
Let $0 < \alpha \leq 1$, $2 \alpha \leq \beta$, and $\lambda > 0$.
Then $t^\alpha E_{\alpha,\beta} (-\lambda t^\alpha)$ is a monotonically increasing function of $t > 0$.
\end{corollary}

\proof
Since $\beta-\alpha \geq \alpha$, then by Lemma \ref{lem:mlf-cm},
$E_{\alpha,\beta-\alpha}(-\lambda t^\alpha)$ is a monotonically decreasing function of $t > 0$.
Let $0 < s < t$, then from \eqref{eq:ML-recurrence} we have
\begin{eqnarray*}
\lambda s^\alpha E_{\alpha,\beta}(-\lambda s^\alpha)
&=&
1/\Gamma(\beta-\alpha) - E_{\alpha,\beta-\alpha}(-\lambda s^\alpha)
\\ &<&
1/\Gamma(\beta-\alpha) - E_{\alpha,\beta-\alpha}(-\lambda t^\alpha)
=
\lambda t^\alpha E_{\alpha,\beta}(-\lambda t^\alpha).  \eproof
\end{eqnarray*}

In the coming sections, we deal with a $2 \times 2$ linear system with a coefficient matrix of the form
\begin{equation}
\label{eq:A-def}
A_{\alpha,\gamma,\mu}(\lambda,s,t):=
\left[ \begin{array}{cc}
t^\alpha E_{\alpha,\mu}(-\lambda t^\alpha) &
t^{\gamma-1}E_{\alpha,\gamma}(-\lambda t^\alpha)
\\
s^\alpha E_{\alpha,\mu}(-\lambda s^\alpha) &
s^{\gamma-1}E_{\alpha,\gamma}(-\lambda s^\alpha)
\end{array} \right].
\end{equation}

In the next lemma, we show that the determinant of this matrix, denoted by $|\dots|$, has a positive lower bound, and consequently is non-singular. This property plays a crucial role for obtaining the coefficients of the series representation of $u$ and $f$ in the forthcoming sections. Also, it provides the necessary bounds for showing the convergence of these series.

\begin{lemma}
\label{lem:matrix}
Let $0 < \alpha \leq \gamma \leq 1$, $\mu \geq 2\alpha$, and $0 < s < t$.
Then, there is a constant $A>0$, independent of $\lambda$, such that
\begin{equation}
\label{eq:Abound}
|A_{\alpha,\gamma,\mu}(\lambda,s,t)| > \frac{A}{\lambda^2}, \qquad \lambda > 0.
\end{equation}
\end{lemma}

\proof
By Lemma \ref{lem:mlf-cm} and Corollary \ref{cor:mlf-increasing},
from \eqref{eq:A-def}, the determinant of the matrix is
$$
\left| A_{\alpha,\mu,\gamma}(\lambda,s,t) \right| =
t^\alpha E_{\alpha,\mu}(-\lambda t^\alpha)
\;
s^{\gamma-1} E_{\alpha,\gamma}(-\lambda s^\alpha)
-
t^{\gamma-1} E_{\alpha,\gamma}(-\lambda t^\alpha)
\;
s^\alpha E_{\alpha,\mu}(-\lambda s^\alpha)>0.
$$
In addition, from \eqref{eq:ML-recurrence} we have
\begin{eqnarray*}
\lim_{\lambda  \rightarrow \infty}
\lambda^2| A_{\alpha,\gamma,\mu}(\lambda,s,t)|
&=&
\frac{1}{\Gamma(\mu-\alpha)}
\frac{s^{\gamma-1-\alpha}}{\Gamma(\gamma-\alpha)}
-
\frac{t^{\gamma-1-\alpha}}{\Gamma(\gamma-\alpha)}
\frac{1}{\Gamma(\mu-\alpha)}.
\\ &=&
\frac{s^{\gamma-1-\alpha} - t^{\gamma-1-\alpha}}{\Gamma(\gamma-\alpha) \;  \Gamma(\mu-\alpha)}
> 0.
\end{eqnarray*}
Therefore, $\lambda^2| A_{\alpha,\gamma,\mu}(\lambda,s,t)| $, as a function of $\lambda$, is bounded below by a positive constant. \eproof

By combining Lemmas \ref{lem:mlf-bound} and \ref{lem:matrix} we obtain the following estimate.

\begin{corollary}
Let $0 < \alpha \leq \gamma \leq 1$, $\mu \geq 2\alpha$, and $0 < s < t$.
Then, there is a constant $B>0$, independent of $\lambda$, such that
\begin{equation}
\label{eq:Ainv-bound}
\left[ A^{-1}_{\alpha,\gamma,\mu}(\lambda,s,t) \right]_{ij} < B \lambda,
\qquad \lambda > 0.
\end{equation}

\end{corollary}

We conclude this section by the following Laplace transform formulas.
\begin{equation}
\label{eq:laplace-general}
\lap \left\{t^{\beta-1} E_{\alpha,\beta}^\rho(\lambda t^\alpha) \right\}
= \frac{s^{\alpha \rho -\beta}}{(s^\alpha-\lambda)^\rho},
\quad|\lambda s^{-\alpha}|<1,
\end{equation}
where $\alpha, \beta,\lambda, \rho \in \C$, $\Re \alpha >0$, $\Re \beta >0$, and $\Re \rho >0$.
See for example (1.9.13) in \cite{kilbas:06b} and (11.8) in \cite{haubold:11}.
From \cite{hilfer:00evolution}, we have the formula
\begin{equation}
\label{eq:laplace-D}
\lap \left\{\hD w(t) \right\} = s^\alpha \lap\{w\} - s^{\alpha-\gamma} I^{1-\gamma} w(0).
\end{equation}

\section{Series Representations}
\label{sec:soln}

Following \cite{furati:14,kirane:11},
the boundary conditions in \eqref{eq:heat-equation} suggest the bi-orthogonal pair of bases $\Phi =  \{\varphi_{0}, \varphi_{1,n}, \varphi_{2,n}\}_{n=1}^\infty$ and $\Psi =  \{\psi_{0}, \psi_{1,n}, \psi_{2,n}\}_{n=1}^\infty$ for the space $L^2(0,1)$
where,
 \begin{equation}
\label{eq:basis1}
\varphi_{1,0}(x) = 2(1-x), \qquad \varphi_{1,n}(x) = 4(1-x)\cos \lambda_n x,
\qquad \varphi_{2,n}(x) = 4\sin \lambda_n x,
\end{equation}
with $\lambda_n = 2\pi  n $, and
\begin{equation}
\label{eq:basis2}
\psi_{1,0}(x) = 1, \qquad \psi_{1,n}(x) = \cos \lambda_n x,
\qquad \psi_{2,n}(x) = x \sin \lambda_n x.
\end{equation}
Although the sequences $\Phi$ and $\Psi$ are not orthogonal, it is proven in \cite{ilin:03} that they both are Riesz bases.  We seek series representations of the solution $u$ and the source term $f$ in the form
\begin{equation}
\label{u-expansion}
u(x,t) =
u_{1,0}(t) \, \varphi_{1,0}(x) + \sum_{k=1}^2\sum_{n=1}^{\infty} u_{k,n}(t) \, \varphi_{k,n}(x),
\end{equation}
\begin{equation}
\label{f-expansion}
f(x)= f_{1,0} \, \varphi_{1,0}(x) + \sum_{k=1}^2\sum_{n=1}^{\infty} f_{k,n} \, \varphi_{k,n}(x).
\end{equation}
Substituting \eqref{u-expansion} and \eqref{f-expansion} into \eqref{eq:heat-equation} yields the following system of fractional differential equations:
\begin{align}
\label{eq:u1}
\hD u_{1,n}(t)  +\lambda_n^2 u_{1,n}(t) &= f_{1,n}, \qquad \qquad n\ge 0,
\\
\label{eq:u2}
\hD u_{2,n}(t) + \lambda_n^2 u_{2,n}(t) - 2 \lambda_n u_{1,n}(t) &= f_{2,n},\qquad \qquad n\ge 1,
\end{align}
where $\lambda_0:=0$.

We solve \eqref{eq:u1} and \eqref{eq:u2} via Laplace transform.
For convenience, we introduce the following notations:
\begin{equation*}
c_{1,0} = I^{1-\gamma} u_{1,0}(0)\quad{\rm and}\quad c_{k,n} = I^{1-\gamma}u_{k,n}(0), \quad k=1,2,\quad n\ge 1.
\end{equation*}
By applying the Laplace transform \eqref{eq:laplace-D} to \eqref{eq:u1}, we obtain
\begin{equation}
\label{eq:U1n}
U_{1,n}(s) =
f_{1,n} \, \frac{1}{s(s^{\alpha} +\lambda_n^2)}  +
c_{1,n} \, \frac{s^{\alpha-\gamma}}{s^\alpha + \lambda_n^2}.
\end{equation}
Similarly, by applying the Laplace transform to \eqref{eq:u2} and then substituting \eqref{eq:U1n}, we obtain
\begin{equation*}
\begin{aligned}
U_{2,n}(s) &=
\frac{f_{2,n}}{s(s^\alpha+\lambda_n^2)} +
\frac{c_{2,n} \, s^{\alpha-\gamma}}{s^\alpha+\lambda_n^2}  +
\frac{2 \lambda_n}{s^\alpha+\lambda_n^2} \, U_{1,n}(s)
\\ &=
\frac{f_{2,n}}{s(s^\alpha+\lambda_n^2)} +
\frac{c_{2,n} \, s^{\alpha - \gamma}}{s^\alpha+\lambda_n^2} +
2 \lambda_n \left[ f_{1,n} \, \frac{1}{s(s^\alpha +\lambda_n^2)^2} +
c_{1,n}  \, \frac{s^{\alpha-\gamma}}{(s^\alpha +\lambda_n^2)^2} \right].
\end{aligned}
\end{equation*}
Hence, from formula \eqref{eq:laplace-general} we have
\begin{equation}
\label{u1n-solution}
u_{1,n}(t) = f_{1,n} \, t^\alpha E_{\alpha,\alpha+1}(-\lambda_n^2 t^\alpha) +
c_{1,n} \,  t^{\gamma-1}E_{\alpha,\gamma}(-\lambda_n^2 t^\alpha),
\qquad n\ge 0,
\end{equation}
and
\begin{equation}
\label{u2n-solution}
u_{2,n}(t) = f_{2,n} \, t^\alpha E_{\alpha,\alpha+1}(-\lambda_n^2 t^\alpha) +
c_{2,n} \, t^{\gamma-1}E_{\alpha,\gamma}(-\lambda_n^2 t^\alpha)
+ d_n(t),\qquad n\ge 1,
\end{equation}
where
\begin{equation}
\label{eq:dn}
d_n(t) = 2 \lambda_n \left[
f_{1,n} \, t^{2\alpha} \, E_{\alpha,2\alpha+1}^2 (-\lambda_n^2 t^\alpha)
+
c_{1,n} \, t^{\alpha+\gamma-1} \,
E_{\alpha,\alpha+\gamma}^2 (-\lambda_n^2 t^\alpha) \right].
\end{equation}

Next, we determine the unknowns
$\{c_{1,0}, c_{k,n}\}$ and $\{f_{1,0}$, $f_{k,n}\}$.
From the two time conditions in \eqref{eq:heat-equation} we have
\begin{equation}
\label{eq:ukn(T)}
u_{1,0} (T_m) = z_{1,0}, \qquad  u_{k,n} (T_m) = z_{k,n},\qquad
u_{1,0} (T) = h_{1,0}, \qquad u_{k,n} (T) = h_{k,n},
\end{equation}
for $n \ge 1$ with $k=1,\,2$, where $\{z_{1,0}, z_{k,n}\}$ and  $\{h_{1,0}, h_{k,n}\}$ denote the Fourier coefficients of the series representations of $z$ and  $h$ in terms of the basis~\eqref{eq:basis1}, respectively. 
That is,
$$
z_{1,0}= \langle z,\psi_{1,0} \rangle , \qquad z_{k,n}= \langle z,\psi_{k,n} \rangle,
\qquad
 h_{1,0} = \langle h, \psi_{1,0} \rangle,\qquad h_{k,n} = \langle h, \psi_{k,n} \rangle,
$$
where $\langle \cdot,\cdot \rangle$ denotes the inner product in $L^2(0,1)$.

Using \eqref{eq:ukn(T)} in \eqref{u1n-solution}, we obtain the linear system
\begin{equation}
\label{eq:system}
A_n
\left[ \begin{array}{c} f_{1,n} \\ c_{1,n} \end{array} \right] =
\left[ \begin{array}{c} h_{1,n} \\ z_{1,n} \end{array} \right],
\qquad n\geq 0,
\end{equation}
where,
\begin{equation}
\label{eq:An}
A_n = A_{\alpha, \gamma, 1 + \alpha}(\lambda_n^2, T_m, T) =
\left[ \begin{array}{cc}
T^\alpha E_{\alpha,\alpha+1}(-\lambda_n^2 T^\alpha) &
T^{\gamma-1}E_{\alpha,\gamma}(-\lambda_n^2 T^\alpha)
\\
T_m^\alpha E_{\alpha, \alpha+1}(-\lambda_n^2 T_m^\alpha) &
T_m^{\gamma-1}E_{\alpha,\gamma}(-\lambda_n^2 T_m^\alpha)
\end{array} \right].
\end{equation}
By Lemma \ref{lem:matrix}, the linear system \eqref{eq:system} is uniquely solvable and
\begin{equation}
\label{eq:f1c1}
\left[ \begin{array}{c} f_{1,n} \\ c_{1,n} \end{array} \right] =
A_n^{-1} \;
\left[ \begin{array}{c} h_{1,n} \\ z_{1,n} \end{array} \right],
\qquad n\geq 0.
\end{equation}
In a similar fashion, by  using \eqref{eq:ukn(T)} in \eqref{u2n-solution}, we observe
\begin{equation}
\label{eq:f2c2}
\left[ \begin{array}{c}
f_{2,n} \\ c_{2,n} \end{array} \right] =
A_n^{-1}
\left[ \begin{array}{c}
h_{2,n} - d_n(T) \\ z_{2,n} - d_n(T_m)
\end{array} \right],
\qquad n\ge 1.
\end{equation}

To determine the coefficients in the series representations
\eqref{u-expansion} and \eqref{f-expansion}, we first perform the calculations in \eqref{eq:f1c1} and \eqref{eq:f2c2} and then substitute in the formulas \eqref{u1n-solution} and \eqref{u2n-solution}.  Therefore, the construction of the series representation of $u$ and $f$ is now completed.

\section{Existence and Uniqueness}
\label{sec:results}

In this section, under some regularity assumptions on the given data functions $h$ and $z$ in problem \eqref{eq:heat-equation}, we show that the series representations of the solution $u$ in \eqref{u-expansion} and of the source $f$ in \eqref{f-expansion} satisfy certain smoothness properties. These smoothness properties will allow us to show the existence and uniqueness of such $u$ and $f$, also  to show that $u$ form a classical solution of \eqref{eq:heat-equation}.


\begin{theorem}
\label{thm:existence}
Let $h,\,z \in C^4[0,1]$ be such that
\begin{equation}
\label{eq:zh}
z'(0) = z'(1), \qquad z(1) = z''(1) = 0, \qquad
h'(0) = h'(1), \qquad h(1) = h''(1) = 0.
\end{equation}
Let $u$ and $f$ be as determined in the previous section.
Then $u(.,t) \in C^2[0, 1]$, $\hD u(x,.) \in C(0,T]$, and $f \in C[0,1]$.
In addition, $u$ and $f$ form the unique classical solution and source of
\eqref{eq:heat-equation}, respectively.
\end{theorem}

\proof
Let $\Omega = (0,T] \times [0,1]$ and $\Omega_\epsilon = [\epsilon,T] \times [0,1] \subset \Omega$.
We show that the series corresponding to $u$, $u_x$, $u_{xx}$, $D^{\alpha,\gamma} u$ are uniformly convergent and represent continuous functions on $\Omega_\epsilon$, for any $\epsilon > 0$.  Also we show that the series representation of $f$ is uniformly convergent in $[0,1]$.  This is shown by bounding all these series by over-harmonic series then applying Weierstrass M-test.
Throughout this proof, $L = L(\alpha,\gamma)$ is some positive (generic) constant.

Through repeated integration by parts, the assumptions in \eqref{eq:zh} yield
$$
\langle z, \psi_{1,n} \rangle = \frac{1}{\lambda_n^4}
\left \{
 z'''(0) - z'''(1) + \langle z^{(4)},\psi_{1,n} \rangle  \right\},
$$
and
$$
\langle z, \psi_{2,n} \rangle = \frac{1}{\lambda_n^4}
\int_0^1 \left[ 4 z'''(x) + xz^{(4)}(x) \right] \sin \lambda_n x \; dx.
$$
Same expressions hold for the function $h$.
Thus, there is a constant $L > 0$ such that
\begin{equation}
\label{eq:zh-bound}
\left| z_{k,n} \right| \leq L/\lambda_n^4, \qquad
\left| h_{k,n} \right| \leq L/\lambda_n^4,
\qquad n\ge 1, \quad  k=1,2.
\end{equation}

Using these bounds and \eqref{eq:Ainv-bound}, it follows from \eqref{eq:f1c1} that
\begin{equation}
\label{eq:c1f1-bound}
|c_{1,n}| \leq L/\lambda_n^2, \qquad |f_{1,n}| \leq L/\lambda_n^2, \qquad n\ge 1.
\end{equation}
Consequently, from \eqref{eq:dn} and \eqref{eq:mlf-bound} we have
$|d_n(T)|, |d_n(T_m)| \leq L/ \lambda_n^5$, $n \geq 1$.
Therefore, it follows from \eqref{eq:f2c2} that
\begin{equation}
\label{eq:c2f2-bound}
|c_{2,n}| \leq L/\lambda_n^2, \qquad |f_{2,n}| \leq L/\lambda_n^2, \qquad n\ge 1.
\end{equation}

Using the bounds \eqref{eq:zh-bound}, \eqref{eq:c1f1-bound}, \eqref{eq:c2f2-bound} and \eqref{eq:mlf-bound}, the formulas \eqref{u1n-solution} and \eqref{u2n-solution} imply that
\begin{equation}
\label{eq:ukn-bound}
t^{1-\gamma+\alpha} |u_{k,n}(t)| \leq L \lambda_n^4, \qquad  n\ge 1, \quad k=1,2, \quad t \in (0,T].
\end{equation}
Furthermore, by inserting this bound in  the \eqref{eq:u1}-\eqref{eq:u2}, we have
\begin{equation}
\label{eq:Dukn-bound}
t^{1-\gamma+\alpha} \, |\hD u_{k,n}(t)| \leq L/\lambda_n^2,
\qquad  n\ge 1, \quad k=1,2,\quad t \in (0,T].
\end{equation}

When $n=0$, it is obvious from \eqref{u1n-solution} that the first term
$u_{1,0}(t) \phi_{1,0}(x)$ of the series \eqref{u-expansion} is continuous
on $\Omega$ and has a continuous $t$ fractional derivative on $\Omega$.
In addition, it has a continuous first and second derivative with respect to $x$ on $\Omega$.

Therefore, the series in \eqref{u-expansion} is uniformly convergent in $\Omega_\epsilon$.
Furthermore, the series obtained through term-by-term fractional differentiation with respect to $t$, and through term-by-term first and second differentiation with respect to $x$ are all uniformly convergent in $\Omega_\epsilon$.  Hence, being represented by uniformly convergent series of continuous functions on $\Omega_\epsilon$, the functions $u$, $\hD u$, $u_x$, and $u_{xx}$ are all continuous on $\Omega$.
Similarly, $f$ is continuous on $[0,1]$.

The next step is to show that $u(x,t)$ satisfies the intermediate and final conditions.  From \eqref{eq:ukn(T)}, the series representation \eqref{u-expansion} at $t = T_m$ yields
$$
u(x,t)|_{t=T_m} =
2z_{1,0}(1-x) +
\sum_{n=1}^{\infty} 4(1-x)z_{1,n}\cos(\lambda_n x) +
\sum_{n=1}^{\infty} 4z_{2,n}\sin(\lambda_n x) =
z(x), \quad x \in [0,1].
$$
Similarly, we have $u(x,T) = h(x)$, $x \in [0,1]$.
Therefore, $u$ form a classical solution.

Finally, the uniqueness follows by observing that when $z=h=0$ in
\eqref{eq:heat-equation}, then from \eqref{eq:f1c1} and \eqref{eq:f2c2}, all the coefficients $c_{kn}$ and $f_{kn}$ are zero.  As a result, from
 \eqref{u1n-solution} and \eqref{u2n-solution}, all the coefficients in the series representations of $u$ and $f$ are identically zero.  \eproof

\section{Analytical and Computational Examples}
\label{sec:example}
To complement the achieved results, an analytical and a numerical example are presented.

\subsection{Example 1 (Linear source and source-free diffusion)}
Consider the problem \eqref{eq:heat-equation} with
$$
z(x) = 2(1-x), \quad \text{and} \quad h(x) = 2c(1-x), \quad c>0.
$$
Then clearly, $z$ and $h$ satisfy the hypothesis of Theorem \ref{thm:existence}, and
$$
z_{1,0} = 1, \quad h_{1,0}  = c, \qquad
z_{k,n} =  h_{k,n} = 0, \qquad n\ge 1, \quad k=1,2.
$$
Thus, it follows from \eqref{eq:f1c1} and  \eqref{eq:dn} that
$$
c_{1,n} = f_{1,n} = d_n  = 0, \qquad n \geq 1,
$$
which imply that
$$
c_{2,n} = f_{2,n} = 0, \qquad n \geq 1.
$$
Accordingly, from \eqref{eq:u1} and \eqref{eq:u2}, for $t \in (0,T]$,
$$
u_{k,n}(t) = 0, \qquad n \geq 1, \quad k=1, 2.
$$
Thus, from \eqref{u-expansion} and \eqref{f-expansion}, the series representations of $f$ and $u$  reduce to
$$
f(x) = 2 f_{1,0}(1-x),
$$
and
$$
u(x,t) = 2 u_{1,0}(t)(1-x) =
2 \left[ \frac{f_{1,0}}{\Gamma(\alpha+1)} \; t^\alpha
+ \frac{c_{1,0}}{\Gamma(\gamma)} \; t^{\gamma-1} \right] (1-x).
$$

Next, from \eqref{eq:f1c1}, to determine  $f_{1,0}$ and $c_{1,0}$ we
calculate $A_0^{-1}$.
From \eqref{eq:An},
$$
A_0 =
\left[ \begin{array}{cc}
T^\alpha E_{\alpha,\alpha+1}(0) &
T^{\gamma-1}E_{\alpha,\gamma}(0)
\\
T_m^\alpha E_{\alpha, \alpha+1}(0) &
T_m^{\gamma-1}E_{\alpha,\gamma}(0)
\end{array} \right]
=
\left[ \begin{array}{cc}
T^\alpha/\Gamma(\alpha+1) &
T^{\gamma-1} / \Gamma(\gamma)
\\
T_m^\alpha /\Gamma(\alpha+1) &
T_m^{\gamma-1} / \Gamma(\gamma)
\end{array} \right].
$$
Thus,
$$
A^{-1}_0
=
\frac{1}{T_m^{\gamma-1} T^\alpha - T^{\gamma-1} T_m^\alpha}
\left[ \begin{array}{cc}
\Gamma(\alpha + 1) \; T_m^{\gamma-1} &
- \Gamma(\alpha + 1) \; T^{\gamma-1}
\\
- \Gamma(\gamma) \; T_m^\alpha &
\Gamma(\gamma) \; T^\alpha
\end{array} \right].
$$
Therefore,
$$
\left[ \begin{array}{c} f_{1,0} \\ c_{1,0} \end{array} \right] =
A^{-1}_0 \left[ \begin{array}{c} h_{1,0} \\ z_{1,0} \end{array} \right]
=
A^{-1}_0 \left[ \begin{array}{c} c \\ 1 \end{array} \right]
=
\frac{1}{T_m^{\gamma-1} T^\alpha - T^{\gamma-1} T_m^\alpha}
\left[ \begin{array}{c}
\Gamma(\alpha + 1) ( c T_m^{\gamma-1} -  T^{\gamma-1} )
\\
\Gamma(\gamma) ( T^\alpha -c T_m^\alpha )
\end{array} \right].
$$

Notice that, when $c = (T_m/T)^{1-\gamma}$, then $f(x)=0$ and thus problem
\eqref{eq:heat-equation} is source-free.
On the other hand, when $c = (T/T_m)^\alpha$, then the problem corresponds to the homogeneous initial condition, $u(x,0)= 0$.

\subsection{Example 2}
 Consider the problem \eqref{eq:heat-equation} with
$$
z(x)= \frac{1}{10} \, \left[1 - (2x-1)^3 \right]
\quad \text{and} \quad
h(x) =  -3x^2(x^2+2)+8x^3+1.
$$
Then, by direct calculations, we can verify that  $z$ and $h$ satisfy the hypothesis of Theorem~\ref{thm:existence}.

We solve the system of equations in \eqref{eq:f1c1} and \eqref{eq:f2c2} to find the Fourier series coefficients and then substitute in the series representations of $u$ and $f$ in \eqref{u-expansion} and \eqref{f-expansion}. We evaluated $u$ and $f$ by truncating the series in \eqref{u-expansion} and \eqref{f-expansion} after 20 terms.

\begin{figure}
\begin{subfigure}{.5\textwidth}
  \centering
  \includegraphics[width=1.06\linewidth]{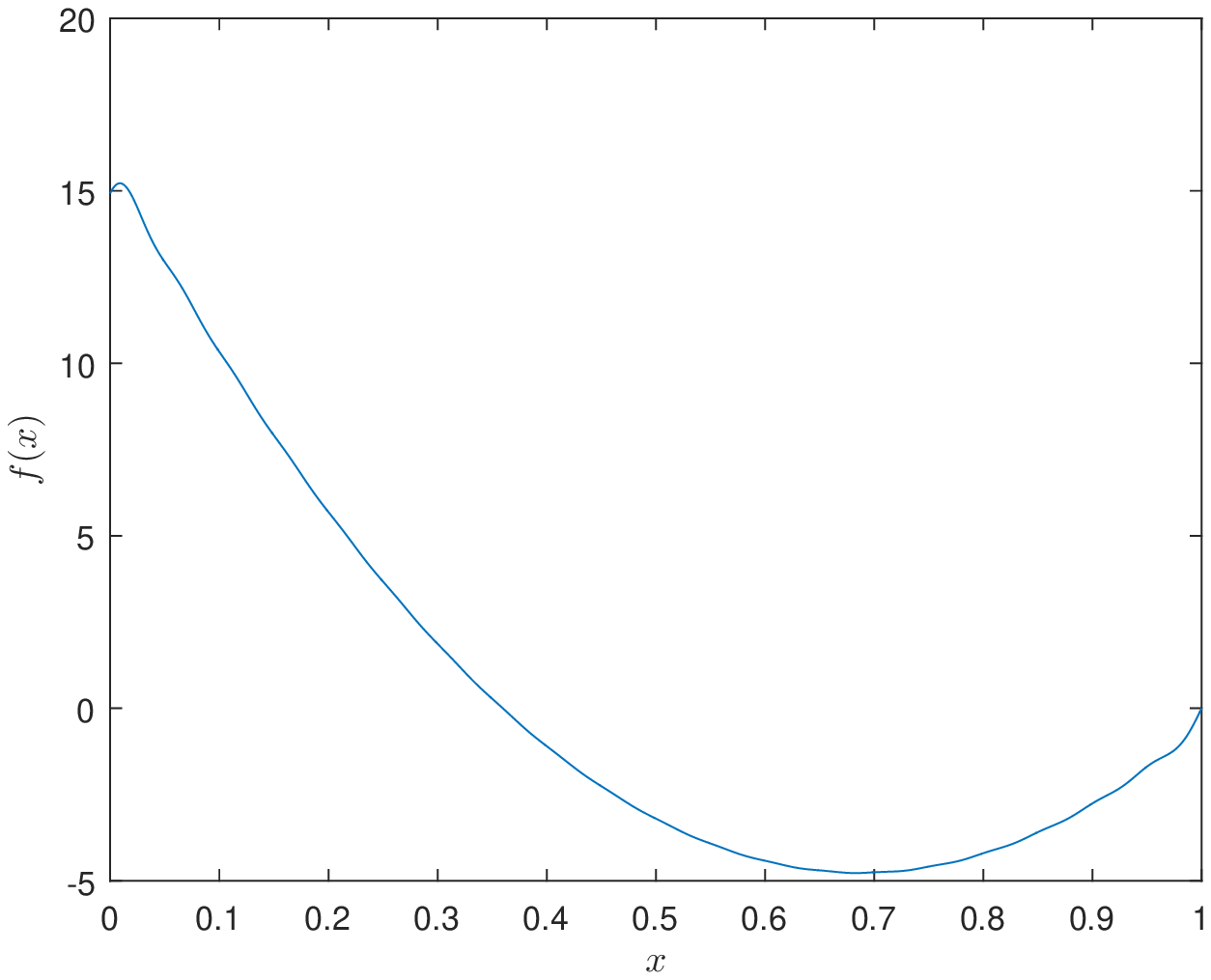}
  \label{fig2:sfig1}
\end{subfigure}%
\begin{subfigure}{.5\textwidth}
  \centering
  \includegraphics[width=1.06\linewidth]{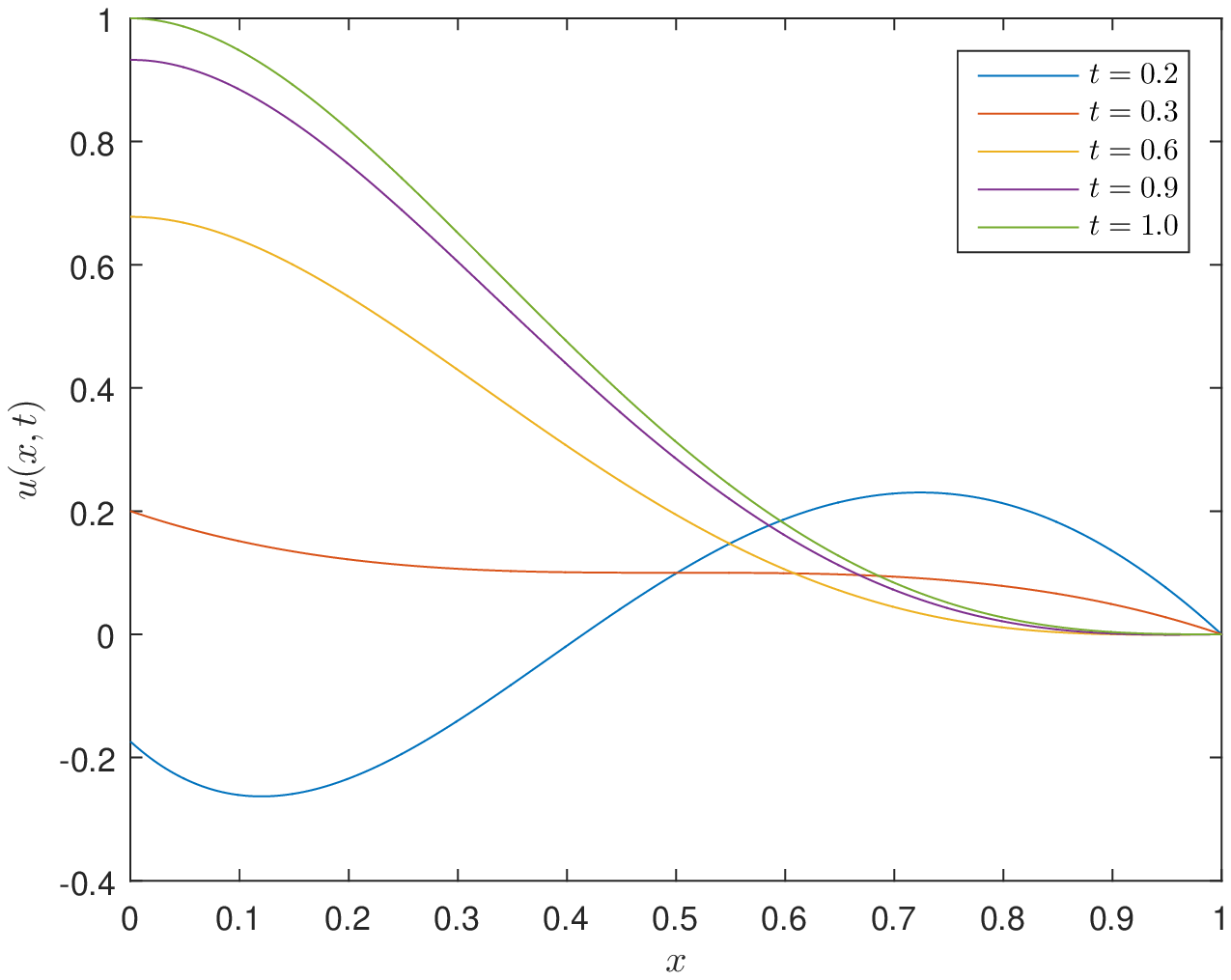}
  \label{fig2:sfig2}
\end{subfigure}
\caption{The solution and source term for Example 2.}
\label{fig2}
\end{figure}

The graph of $u$ at different times and the graph of $f$ are shown in Figure \ref{fig2} for
$T=1$, $T_m= 0.3$, and $\alpha=\gamma=0.5$.
The graph of $u$ at $t =0.2$ represents the solution prior to the measurement at $T_m = 0.3$.

\section*{Acknowledgment}
The authors would like to acknowledge the support provided by the Deanship of Scientific Research at King Fahd University of Petroleum \& Minerals under Research Grant FT151003.

\bibliographystyle{abbrv}
\bibliography{furati,fractional}

\end{document}